\documentclass[11pt]{article}

\usepackage{amsmath,amscd,amssymb,amsthm}
\usepackage{latexsym,enumerate,graphicx,geometry,hyperref}
\usepackage[all]{xy}

\newtheorem{thm_}{Theorem}[section]
\newtheorem{lemma_}[thm_]{Lemma}
\newtheorem{prop_}[thm_]{Proposition}
\newtheorem{def_}[thm_]{Definition}
\newtheorem{rk_}[thm_]{Remarks}
\newtheorem{cor_}[thm_]{Corollary}
\newcommand{\thm}[1]{\begin{thm_}#1\end{thm_}}
\newcommand{\lemma}[1]{\begin{lemma_}#1\end{lemma_}}
\newcommand{\prop}[1]{\begin{prop_}#1\end{prop_}}
\newcommand{\defi}[1]{\begin{def_}#1\end{def_}}
\newcommand{\rk}[1]{\begin{rk_}#1\end{rk_}}

\newcommand{\pf}[1]{\begin{proof}#1\end{proof}}
\newcommand{\fracl}[3]{\genfrac{(}{)}{}{}{#1}{#2}_{#3}}
\newcommand{\fracp}[2]{\fracl{#1}{#2}{p}}
\newcommand{\ord}[2]{\text{ord}_{#1}({#2})}
\newcommand{\tr}[3]{\text{T}_{{#1}/{#2}}({#3})}
\newcommand{\trn}[2]{\text{T}_{{#1}/{#2}}}
\newcommand{\nm}[3]{\text{N}_{{#1}/{#2}}({#3})}
\newcommand{\nmn}[2]{\text{N}_{{#1}/{#2}}}
\renewcommand{\th}{^{\text{th}}}
\newcommand{\card}[1]{\#({#1})}
\newcommand{\cnm}{\mathcal N}
\newcommand{\s}{\sigma}
\newcommand{\Z}{\mathbb Z}
\newcommand{\QQ}{\mathbb Q}

\newcommand{\R}{\mathbb R}
\renewcommand{\O}{\mathcal O}
\newcommand{\p}{\mathfrak p}
\newcommand{\q}{\mathfrak q}
\renewcommand{\P}{\mathfrak P}
\newcommand{\Q}{\mathfrak Q}
\newcommand{\M}{\mathfrak M}
\renewcommand{\L}{\mathfrak L}
\newcommand{\A}{\mathfrak A}
\newcommand{\B}{\mathfrak B}
\newcommand{\eq}[1]{\begin{equation}#1\end{equation}}
\newcommand{\eqn}[1]{\begin{equation*}#1\end{equation*}}

\newcommand{\gan}[1]{\begin{gather*}#1\end{gather*}}

\newcommand{\aln}[1]{\begin{align*}#1\end{align*}}
\renewcommand{\sp}[1]{\begin{split}#1\end{split}}

\begin{document}
\title{ Primality Test for Numbers of the Form $Ap^n+w_n$}
\author{
Yingpu Deng and Chang Lv\\\\
Key Laboratory of Mathematics Mechanization,\\
NCMIS, Academy of Mathematics and Systems Science,\\
Chinese Academy of Sciences, Beijing 100190, P.R. China\\
Email: \nolinkurl{{dengyp,lvchang}@amss.ac.cn}
}
\date{}
\maketitle

\begin{abstract}
    We propose an algorithm determining the primality of numbers $M=Ap^n+w_n$ where $w_n^{p-1}\equiv1\pmod{p^n}$ and $A<p^n$ and give example when $p=7$. $p\th$ reciprocity law is involved. The algorithm runs in polynomial time in $\log_2(M)$ for fixed $p$ and $A$.
\end{abstract}

\section{Introduction}
    In 1983, Adleman, Pomerance and Rumely \cite{APR} gave a general deterministic primality test which is still very practical now. This test was simplified later by Cohen and Lenstra \cite{CoLen} and now is called APRCL test. However this is not running in polynomial time. Agrawal, Kayal and Saxena \cite{AKS} discovered a deterministic polynomial time algorithm for general primality tests in 2004. But it is difficult to make use of it in practice. So finding more efficient algorithms for specific families of numbers makes a lot of sense. Primality tests for numbers of the form $Ap^n\pm1$ with $p$ prime, have been noticed since Lucas \cite{Lucas} and Lehmer \cite{Mersenne} gave the celebrated Lucas-Lehmer primality test for Mersenne numbers, using properties of the Lucas sequences. Here, we recall this famous primality test:

    \textbf{Lucas-Lehmer test.}\quad Let $M_p=2^p-1$ be Mersenne number, where $p$ is an odd prime. Define a sequence $\{u_k\}$ as follows: $u_0=4$ and $u_k=u_{k-1}^2-2$ for $k\geq1$. Then $M_p$ is a prime if and only if $u_{p-2}\equiv0\pmod{M_p}$.

    H.C. Williams and collaborators extended this method to $p=3,5,7$ and even general $p$ and gave many concrete algorithms, see \cite{Williams3,Williams3-,Williams35,Williams57,Williamstri,Stein}. A comprehensive treatise on this method can be found in the book by Williams \cite{Williams}. Generally speaking, this method is rather complicated.

    Another classical line for primality test is Proth theorem.

    \textbf{Proth theorem.}\quad Let $N=h\cdot2^n+1$ with $h$ odd and $h<2^n$. Suppose $a$ is an integer with the Jacobi symbol $\left(\frac{a}{N}\right)=-1$. Then $N$ is a prime if and only if $a^{(N-1)/2}\equiv-1\pmod{N}$.

    Recently, primality tests for $Ap^n\pm1$ based on Proth's theorem and higher reciprocity law rather than Lucas sequences have been developed for small primes $p$. The first paper of this kind is A. Guthmann's \cite{Guthmann} using cubic reciprocity to deal with the primality of $A3^n+1$. Since then P. Berrizbeitia and collaborators continued this research line and presented primality tests for $p=2$, $3$, $5$ and even general $p$, see \cite{Biquadratic,A3n,A5n,A5n2,Amn+wn}. In \cite{Amn+wn}, they gave a generalization of Proth theorem in cyclotomic fields, so one needs to do computation in non-rational number field for their primality tests. The reason is that they didn't give the sequence form of their primality tests, likewise the Lucas-Lehmer test for Mersenne numbers. Since we hope to do computation in rational number field, it is our desire to do primality tests involving only computation in rational number field.

    In this paper we present an efficient criteria for general $p$ with explicit sequence form, determining the primality of $Ap^n+w_n$ where $w_n^{p-1}\equiv1\pmod{p^n}$ and $A<p^n$. For a general $p$, our test makes use of $\frac{p-1}{2}$ many sequences and we give the explicit recursive formulas for these sequences. Once the seeds for these sequences are given, our tests can determine the primality of numbers of form $Ap^n+w_n$, and our tests involve only computation in rational number field. We also give the concrete example for $p=7$.

    To begin with, let us see the following proposition first:
\prop{\label{TrghoutProp}
Let $M$ and $n$ be positive integers with $M>1$ and $p$ an odd prime, then the following two conditions are equivalent:
\begin{enumerate}[(i)]
\item $p^n\mid M^{p-1}-1$ and $p^n>\sqrt M$.
\item $M=Ap^n+w_n$ where $w_n$ and $A$ are integers such that $0<w_n<p^n$ and that $w_n^{p-1}\equiv1\pmod{p^n}$, and $0\le A<p^n$.
\end{enumerate}
}
\pf{
We may write $M=Ap^n+w_n$ with $A\geq0$ and $0\le w_n<p^n$. Since $M\equiv w_n\pmod{p^n}$, so $p^n\mid M^{p-1}-1$ is equivalent to $w_n^{p-1}\equiv1\pmod{p^n}$ and $0<w_n<p^n$. Now $p^{2n}> M$ implies $p^{2n}>Ap^n+w_n>Ap^n$, so $A<p^n$. For the converse, $A<p^n$ means $A\le p^n-1$ so $M=Ap^n+w_n<(p^n-1)p^n+p^n=p^{2n}$.
}

By Hensel's lemma, it is easy to find the $p-1$ many values of $w_n$. Thus, given $w_n$ satisfying $w_n^{p-1}\equiv1\pmod{p^n}$, there is a unique
$x\pmod{p}$ such that $(w_n+xp^n)^{p-1}\equiv1\pmod{p^{n+1}}$ and we obtain $w_{n+1}=w_n+xp^n$. In particular, we have $w_n=a+kp$ with $0<a<p$ and $k\geq0$ such that $a+kp<p^n$.

    Throughout this paper we suppose that $M$ satisfies anyone of the conditions in Proposition \ref{TrghoutProp} and $A>0$, i.e. $M>p^n$. We will deal with the equation $w_n^{p-1}\equiv1\pmod{p^n}$ by another method in Section \ref{SecComp.Solve} later.

This paper is organized as follows. In Section \ref{SecPre} we give without proof the facts we need, mainly the properties of power residue symbol. In Section \ref{SecMain} we state and prove our main result. In Section \ref{SecComp}, we give some computational considerations of our method.
In Section \ref{Sec.p=7} we give the concrete example for $p=7$. In Section \ref{SecImple} we give the implementation and computational results for $p=3,5$ and $7$.

\section{Preliminaries}\label{SecPre}
    In this section, we fix some notations and state some well-known facts we need later.

    Let $\zeta=\zeta_p=e^{2\pi\sqrt{-1}/p}$ be a primitive $p\th$ root of unity, then we are working in the corresponding cyclotomic field $L=\QQ(\zeta)$ and denote the ring of algebraic integers in $L$ by $\O_L$. Let $K=L\cap\R=\QQ(\zeta+\zeta^{-1})$ be the maximal real subfield of $L$. We know $[L:\QQ]=p-1$ and $K/\QQ$ is also a Galois extension of degree $r:=\frac{p-1}{2}$. Let $G=Gal(L/\QQ)\cong(\Z/p\Z)^*$ so for every integer $c$ with $p\nmid c$ denote by $\s_c$ the element of $G$ that sends $\zeta$ to $\zeta^c$. For $\delta$ in the group ring $\Z[G]$ and $\alpha$ in $L$ with $\alpha\neq0$ we denote by $\alpha^\delta$ the action of $\delta$ on $\alpha$, that is,
\eqn{
\alpha^\delta:=\prod_{\s\in G}\s(\alpha)^{k_{\s}}, \text{ if }\delta=\sum_{\s\in G}k_{\s}\s\text{ where }k_{\s}\in\Z.
}
We know that $Gal(L/K)=\{\s_1,\s_{-1}\}$ and we also write $\s_1=1$ in $\Z[G]$.
    Now we briefly introduce without proof what we will use later. See \cite[Chapter 14]{GTM84} for details. Let $\P$ be a prime ideal in $\O_L$ not divide $p$, then for $\alpha\in\O_L-\P$ there is a unique $p\th$ root of unity $\zeta^j$ with $j\in\Z$ such that
\eqn{
\alpha^{\frac{\cnm\P-1}{p}}\equiv\zeta^j\pmod{\P},
}
where $\cnm\P=\card{\O_L/\P}$ is the absolute norm of $\P$. We define this $\zeta^j$ to be $\fracp{\alpha}{\P}$, called {\em$p\th$ power residue symbol} and also adopt $\fracp{\alpha}{\P}=0$ when $\alpha\in\P$. Here are some properties of it:
\prop{\label{pResSym1}
$\notag$
\begin{enumerate}[(i)]
\item For $\alpha\in\O_L-\P,\fracp{\alpha}{\P}=1$ iff $x^p\equiv\alpha\pmod{\P}$ has a solution in $\O_L$.
\item For all $\alpha\in\O_L,\,\alpha^{\frac{\cnm\P-1}{p}}\equiv\fracp{\alpha}{\P}\pmod{\P}$.
\item $\fracp{\alpha\beta}{\P}=\fracp{\alpha}{\P}\fracp{\beta}{\P}$.
\item If $\alpha\equiv\beta\pmod{\P}$ then $\fracp{\alpha}{\P}=\fracp{\beta}{\P}$.
\end{enumerate}
}
\defi{
Suppose $\A$ is an ideal of $\O_L$ that is prime to $p$. Let $\A=\prod\P$ be the prime decomposition of $\A$. Define $\fracp{\alpha}{\A}=\prod\fracp{\alpha}{\P}$. If $\beta\in\O_L$ and is prime to $p$ define $\fracp{\alpha}{\beta}=\fracp{\alpha}{(\beta)}$.
}
\prop{\label{pResSym2}
$\notag$
\begin{enumerate}[(i)]
\item $\fracp{\alpha\beta}{\A}=\fracp{\alpha}{\A}\fracp{\beta}{\A}$.
\item $\fracp{\alpha}{\A\B}=\fracp{\alpha}{\A}\fracp{\alpha}{\B}$.
\item Let $\s\in G$ then $\fracp{\alpha}{\A}^{\s}=\fracp{\s(\alpha)}{\s(\A)}$.
\end{enumerate}
}
\defi{[Primary element]\label{Primary}
A nonzero element $\alpha\in\O_L$ is called primary if it is not a unit and is prime to $p$ and congruent to a rational integer modulo $(1-\zeta)^2$.
}
The following lemma shows that primary elements are plentiful.
\lemma{
Suppose $\alpha\in\O_L$ and is prime to $p$. There is an integer $c\in\Z$, unique modulo $p$, such that $\zeta^c\alpha$ is primary.
}
Now we state the
\thm{[The Eisenstein Reciprocity Law]\label{Eisenstein}
Let $p$ be an odd prime, $a\in\Z$ prime to $p$, and $\alpha\in\O_L$ a primary element prime to $a$. Then
\eqn{
\fracp{\alpha}{a}=\fracp{a}{\alpha}.
}
}

\section{The Main Theorem}\label{SecMain}
    For $M$ as described before we find a $\pi\in\O_L$ prime to $M$ such that $\fracp{\pi}{M}\neq1$. We will discuss how to find such $\pi$ later in Section \ref{SecComp.Find}.

    By assumption $M$ is coprime with $p$, so let $f=\ord{p}{M}$ be the order of $M$ modulo $p$. Since $f\mid p-1$ and $G$ is cyclic, let $H$ be the unique subgroup of $G$ of order $f$. We know that $H=\{\s_{M^j}\mid 0\le j\le f-1\}$. Denote by $\Phi_f(x)$ the $f\th$ cyclotomic polynomial over $\QQ$, then we have
\prop{\label{pnDivPhi}
$p^n\mid\Phi_f(M)$.
}
\pf{
If $f=1$, since $w_n=1$, the result is obvious. Suppose $f>1$ then $f=\ord{p}{M}$ implies $p\mid M^f-1$. From $M^a-1=\prod_{d\mid a}\Phi_d(M)$ for any integer $a>0$, it follows that $p\mid\Phi_f(M)$ but $p\nmid\Phi_d(M)$ for any $d\mid f$ with $d<f$. However $p^n\mid M^{p-1}-1=(M^f-1)((M^f)^{\frac{p-1}{f}}+\dots+1)$ but $p$ does not divide the second factor, otherwise $0\equiv(M^f)^{\frac{p-1}{f}}+\dots+1)\equiv\frac{p-1}{f}+1\pmod p$, but $0<\frac{p-1}{f}+1<p$, which yields a contradiction. Hence $p^n\mid M^f-1=\prod_{d\mid f}\Phi_d(M)$. By the above argument, we have $p^n\mid\Phi_f(M)$.
}

    Now we distinguish two cases:
\begin{enumerate}[(i)]
\item $f$ is odd. Then we know $\s_{-1}\not\in H$ and $g:=\card{G/H}=\frac{p-1}{f}$ is even. Since $\langle\s_{-1}\rangle H$ is a subgroup of $G$ of order $2f$, it follows that there exist a set of integers $S=\{i_1,i_2,\dots,i_{\frac{g}{2}}\}$ such that $\{\s_{ i_s}\mid s=1,2,\dots,\frac{g}{2}\}$ is a set of representatives of $G/\langle\s_{-1}\rangle H$, hence $\{\s_{\pm i_s}\mid s=1,2,\dots,\frac{g}{2}\}$ is a set of representatives of $G/H$. We set
\eq{\label{gamma1}
\gamma=\sum_{i\in S}i(\s_i^{-1}-\s_{-i}^{-1})\in\Z[G].
}
\item $f$ is even. Let $T=\{j_1,j_2,\dots,j_g\}$ be any set of integers such that $\s_{j_t},\,t=1,2,\dots,g$ is a set of representatives of $G/H$ and set
\eq{\label{gamma2}
\gamma=\sum_{j\in T}j\s_j^{-1}\in\Z[G].
}
\end{enumerate}

    In all cases let
\eq{\label{tau}
\tau=\pi^{\gamma\prod_{d\mid f,d<f}\Phi_d(\s_M)\frac{\Phi_f(M)}{p^n}}
}
which is in $\Z[G]$ by Proposition \ref{pnDivPhi}.

    Let
\eq{\label{tk}
t_k=\tr{L}{K}{\tau^{p^k}}=\tau^{p^k}+\bar\tau^{p^k}\quad\in K
}
where the bar denotes $\s_{-1}$ which acts the same as the complex conjugation and $\trn{L}{K}$ is the trace map from $L$ to $K$. We now turn to a
\prop{
\eq{\label{Tkm}
T_k^{(m)}:=\sum_{1\le i_1<\dots<i_m\le r}\prod_{j=1}^m t_k^{\s_{i_j}},\qquad m=1,\dots,r
}
are all in $\QQ$ for each $k\in\Z$ and $k\geq0$.
}
\pf{
It is easy to see that $Gal(K/\QQ)=\{\s_i|_K\mid i=1,\ldots,r\}$. Now by the symmetric properties of $T_k^{(m)}$ in the definition, all $T_k^{(m)}$ are fixed by each element in $Gal(K/\QQ)$, hence they are in $\QQ$.
}
\rk{
In particular $T_k^{(1)}=\tr{K}{\QQ}{t_k}$ and $T_k^{(r)}=N_{K/\QQ}(t_k)$, where $\trn{K}{\QQ}$ and $\nmn{K}{\QQ}$ are the trace and the norm map from $K$ to $\QQ$, respectively.
}
Recall that $\pi$ is prime to $M$, so $\tau$ could be viewed as in $(\O_L)_{(M)}$ which is the localization of $\O_L$ by $\O_L-\cup_{\P\mid M}\P$. By $a\equiv b\pmod M$ where $a,b\in(\O_L)_{(M)}$, we mean that $a$ and $b$ have the same image under the canonical epimorphism $(\O_L)_{(M)}\longrightarrow\O_L/M\O_L$. Since
\eqn{
(\O_L)_{(M)}\cap K=\bigcap_{\P\mid M}(\O_L)_{\P}\cap K=\bigcap_{\p\mid M}\bigcap_{\P\mid\p}((\O_L)_{\P}\cap K)=\bigcap_{\p\mid M}(\O_K)_{\p}=(\O_K)_{(M)}
}
we may regard $t_k$ as in $(\O_K)_{(M)}$ and similarly $T_k^{(m)}$ as in $\Z_{(M)}=(\O_K)_{(M)}\cap\QQ$. We need another
\prop{
Let $\xi$ be a primitive $p\th$ root of unity and let
\eq{\label{u}
u_{\xi}=\tr{L}{K}{\xi}=\xi+\bar\xi\in K
}
then
\eq{\label{Um}
U_{\xi}^{(m)}:=\sum_{1\le i_1<\dots<i_m\le r}\prod_{j=1}^m u_{\xi}^{\s_{i_j}}\in\QQ,\qquad m=1,\dots,r
}
are independent of the choice of $\xi$, denoted as $U^{(m)}$.
}
\pf{
Let $\xi'$ be another primitive $p\th$ root of unity then $\xi'=\xi^{\s}$ for some $\s\in G$. Since $G$ is abelian
\aln{
U_{\xi}^{(m)}&=(U_{\xi}^{(m)})^{\s}=\sum_{1\le i_1<\dots<i_m\le r}\prod_{j=1}^m (u_{\xi}^{\s})^{\s_{i_j}}\\
&=\sum_{1\le i_1<\dots<i_m\le r}\prod_{j=1}^m\tr{L}{K}{\xi^{\s}}^{\s_{i_j}}=\sum_{1\le i_1<\dots<i_m\le r}\prod_{j=1}^m\tr{L}{K}{\xi'}^{\s_{i_j}}=U_{\xi'}^{(m)}.
}
}
We will give some concrete values for $U^{(m)}$ for small $p$ in Section \ref{SecComp.Um} later.

    Now we can state our main
\thm{\label{MainThm}
Let $M$, $\tau$, $T_k^{(m)}$ and $U^{(m)}$ be as before. Suppose further that $M$ is not divisible by any of the solutions of $x^{p-1}\equiv1\pmod{p^n}$ with $1<x<p^n$. Then the following statements are equivalent:
\begin{enumerate}[(i)]
\item $M$ is prime.
\item There exists a primitive $p\th$ root of unity $\xi$ such that
\eq{\label{MainThm.C2}
\tau^{p^{n-1}}\equiv\xi\pmod M.
}
\item
\eq{\label{MainThm.C3}
T_{n-1}^{(m)}\equiv U^{(m)}\pmod M,\qquad \text{for each }m=1,\dots,r.
}
\end{enumerate}
}
\rk{
If $\pi$ is prime to $M$ but $\pi$ is a $p\th$ power modulo $M$, then through the proof of the theorem later, we still obtain a sufficient condition for primality of $M$:
\prop{
Let $M$, $\tau$, $T_k^{(m)}$ and $U^{(m)}$ be as before. And further suppose $s$ is a positive integer with $p^s>\sqrt{M}$ and $M$ is not divisible by any of the solutions of $x^{p-1}\equiv1\pmod{p^s}$ with $1<x<p^s$. Then the following two statements are equivalent:
\begin{enumerate}[(a)]
\item There exists a primitive $p\th$ root of unity $\xi$ such that
\eqn{
\tau^{p^{s-1}}\equiv\xi\pmod M.
}
\item
\eqn{
T_{s-1}^{(m)}\equiv U^{(m)}\pmod M,\qquad \text{for each }m=1,\dots,r.
}
\end{enumerate}
And at this time, $M$ is a prime.
}
}

    Before we prove the theorem let us see how to deduce a set of recursive formulas for $T_{k}^{(m)}$. Define $F(z_1,z_2)\in\Z[z_1,z_2]$ to be the uniquely determined polynomial such that $x^p+y^p=F(x+y,xy)$, i.e. the representation of the symmetric polynomial $x^p+y^p$ in the elementary symmetric ones. Similarly since
\eqn{
g_m(x_1,\dots,x_r)=\sum_{1\le i_1<\dots<i_m\le r}\prod_{j=1}^mF(x_{i_j},1),\qquad m=1,\dots,r
}
are clearly symmetric polynomials w.r.t $x_1,\dots,x_r$, so we can uniquely define $G_m(z_1,\dots,z_r)$ $\in\Z[z_1,\dots,z_r]$ such that $g_m(x_1,\dots,x_r)=G_m(x_1+\dots+x_r,\dots,x_1\cdots x_r)$. We will give the expressions of $F(z_1,z_2)$ and $G_m(z_1,\dots,x_r)$ in Section \ref{SecComp.F.G} for small $p$ later. Then we have the following
\prop{\label{PropRecurT}
For every $k\ge0$
\eq{
T_{k+1}^{(m)}=G_m(T_k^{(1)},T_k^{(2)},\dots,T_k^{(r)}),\qquad m=1,2,\dots,r.
}
}
The key point is to prove the
\lemma{\label{Neq1}
$N_{L/K}(\tau)=\tau\bar\tau=\tau^{1+\s_{-1}}=1$.
}
\pf{
By the definition \eqref{tau} of $\tau$, it suffice to prove that either $\gamma$ or $\prod_{d\mid f,d<f}\Phi_d(\s_M)$ is annihilated by $1+\s_{-1}$. Clearly $(\s_i^{-1}-\s_{-i}^{-1})(1+\s_{-1})=\s_i^{-1}(1-\s_{-1})(1+\s_{-1})=0$ so if $f$ is odd then by \eqref{gamma1} we have $(1+\s_{-1})\gamma=\sum_{i\in S}i(\s_i^{-1}-\s_{-i}^{-1})(1+\s_{-1})=0$. If $f$ is even then $\s_M^{\frac{f}{2}}-1$ appears in $\prod_{d\mid f,d<f}\Phi_d(\s_M)$. Since $\s_M$ has order $f$, $\s_M^{\frac{f}{2}}$ has order $2$ so $\s_M^{\frac{f}{2}}=\s_{-1}$, which implies $(\s_M^{\frac{f}{2}}-1)(1+\s_{-1})=0$. It follows that $\prod_{d\mid f,d<f}\Phi_d(\s_M)(1+\s_{-1})=0$. This is what we have asserted.
}
\pf{[Proof (of the Proposition \ref{PropRecurT}). ]
Since by definition \eqref{tk}, $t_k=\tau^{p^k}+\bar\tau^{p^k}$, and by Lemma \ref{Neq1} we have $\tau^{p^k}\bar\tau^{p^k}=N_{L/K}(\tau^{p^k})=1$. Thus $t_{k+1}=(\tau^{p^k})^p+(\bar\tau^{p^k})^p=F(\tau^{p^k}+\bar\tau^{p^k},\tau^{p^k}\bar\tau^{p^k})=F(t_k,1)$. Now for each $m=1,\dots,r$, by definition \eqref{Tkm}
\aln{
T_{k+1}^{(m)}&=\sum_{1\le i_1<\dots<i_m\le r}\prod_{j=1}^m t_{k+1}^{\s_{i_j}}=\sum_{1\le i_1<\dots<i_m\le r}\prod_{j=1}^m F(t_k,1)^{\s_{i_j}}=\sum_{1\le i_1<\dots<i_m\le r}\prod_{j=1}^m F(t_k^{\s_{i_j}},1)\\
&=g_m(t_k^{\s_1},\cdots,t_k^{\s_r})=G_m(\sum_{i=1}^r t_k^{\s_i},\dots,\prod_{i=1}^r t_k^{\s_i})=G_m(T_k^{(1)},T_k^{(2)},\dots,T_k^{(r)}).
}
}
    Now we prove the main theorem.
\pf{[Proof (of the main theorem \ref{MainThm}). ]
(i) $\Longrightarrow$ (ii). Since $\fracp{\pi}{M}\neq1$, there exists a primitive $p\th$ root of unity $\xi$ such that $\fracp{\pi}{M}=\xi$. Now $M$ is a rational prime. Let $\M$ be a prime ideal of $\O_L$ lying over $M$, so we know that $f=\ord{p}{M}=f(\M|M)$ is the relative degree since $M$ is unramified in $L$. Recall $H$ is the unique subgroup of $G$ having order $f$, so $H$ is the decomposition group of $M$ in $L$.
\begin{enumerate}[(a)]
\item Suppose first that $f$ is odd. Recall the set of integers $S=\{i_1,\dots,i_{\frac{g}{2}}\}$ is chosen such that $\s_{\pm i_s},\,s=1,\dots,\frac{g}{2}$ is a set of representatives of $G/H$. So we have the decomposition
\eq{\label{Msplits}
M\O_L=\prod_{i\in S}\s_i(\M)\s_{-i}(\M).
}
Then by Proposition \ref{pResSym2}
\eqn{
\sp{
\xi&=\fracp{\pi}{M}=\prod_{i\in S}\fracp{\pi}{\s_i(\M)}\fracp{\pi}{\s_{-i}(\M)}=\prod_{i\in S}\fracp{\pi^{\s_i^{-1}}}{\M}^{\s_i}\fracp{\pi^{\s_{-i}^{-1}}}{\M}^{\s_{-i}}\\
&=\prod_{i\in S}\fracp{\pi^{i\s_i^{-1}}}{\M}\fracp{\pi^{-i\s_{-i}^{-1}}}{\M}=\fracp{\pi^{\gamma}}{\M}
}
}
where we obtain the last equality by definition \eqref{gamma1} with $\gamma=\sum_{i\in S}i(\s_i^{-1}-\s_{-i}^{-1})$. Therefore
\eqn{
\xi=\fracp{\pi^{\gamma}}{\M}\equiv\pi^{\gamma\frac{\cnm\M-1}{p}}=\pi^{\gamma\frac{M^f-1}{p}}=\pi^{\gamma\frac{M^f-1}{\Phi_f(M)}\frac{\Phi_f(M)}{p^n}p^{n-1}}\pmod{\M}.
}
Keep in mind that $\s_M$ is the Frobenius automorphism of $\O_L/\M$ then $\alpha^{\s_M}\equiv\alpha^M\pmod{\M}$ for all $\alpha\in\O_L$. This makes the last equation into
\eqn{
\xi\equiv\pi^{\gamma\prod_{d\mid f,d<f}\Phi_d(\s_M)\frac{\Phi_f(M)}{p^n}p^{n-1}}=\tau^{p^{n-1}}\pmod{\M}\qquad\text{(see the definition of $\tau$ \eqref{tau})}.
}
Now in \eqref{Msplits} we are free to replace $\M$ by $\s(\M)$ for all $\s\in G$ and with the same argument we obtain symmetrically $\xi\equiv\tau^{p^{n-1}}\pmod{\M}$ for all $\M\mid M$. It follows that $\xi\equiv\tau^{p^{n-1}}\pmod M$ since $M$ is unramified in $L$ and we have the composition
\eqn{\begin{CD}
(\O_L)_{(M)}    @>>>    \O_L/M\O_L    @>\sim>>    \bigoplus_{\M\mid M}\O_L/\M.
\end{CD}}

\item If $f$ is even, we have defined before that $T=\{j_1,j_2,\dots,j_g\}$ is a set of integers such that $\s_{j_t},\,t=1,2,\dots,g$ is a set of representatives of $G/H$. Thus
\gan{
M\O_L=\prod_{j\in T}\s_j(\M),\\
\xi=\fracp{\pi}{M}=\prod_{j\in T}\fracp{\pi}{\s_j(\M)}=\prod_{j\in T}\fracp{\pi^{\s_j^{-1}}}{\M}^{\s_j}=\prod_{j\in T}\fracp{\pi^{j\s_j^{-1}}}{\M}=\fracp{\pi^{\gamma}}{\M}.
}
This time $\gamma=\sum_{j\in T}j\s_j^{-1}$ (see definition \eqref{gamma2}). The sequent argument is the same as in (a).
\end{enumerate}

(ii) $\Longrightarrow$ (iii). By \eqref{MainThm.C2} $\tau^{p^{n-1}}\equiv\xi\pmod M$. We also have $\bar\tau^{p^{n-1}}\equiv\bar\xi\pmod M$ so we take trace from $L$ to $K$ to obtain $t_{n-1}=\tr{L}{K}{\tau^{p^{n-1}}}=\tau^{p^{n-1}}+\bar\tau^{p^{n-1}}\equiv\xi+\bar\xi=\tr{L}{K}{\xi}=u_{\xi}\pmod M$. We remark that this congruence should be viewed in $\O_K/M\O_K$, the correctness seen by the commutative diagram
\eqn{\xymatrix{
(\O_K)_{(M)}    \ar@{->>}[d]\ar@{^{(}->}[r]    &(\O_L)_{(M)}    \ar@{->>}[d]\\
\O_K/M\O_K    \ar@{^{(}->}[r]    &\O_L/M\O_L.
}}
Next,
\eqn{
t_{n-1}^{\s_i}\equiv u_{\xi}^{\s_i}\pmod M\quad,\quad i=1,\dots,r.
}
It follows by the definition \eqref{Tkm} and \eqref{Um} that
\eqn{
T_{n-1}^{(m)}=\sum_{1\le i_1<\dots<i_m\le r}\prod_{j=1}^m t_k^{\s_{i_j}}\equiv\sum_{1\le i_1<\dots<i_m\le r}\prod_{j=1}^m u_{\xi}^{\s_{i_j}}=U^{(m)}\pmod{M},\qquad m=1,\dots,r
}
where the congruence is viewed in $\Z/M\Z$. Here we also use the similar commutative diagram
\eqn{\xymatrix{
\Z_{(M)}    \ar@{->>}[d]\ar@{^{(}->}[r]    &(\O_K)_{(M)}    \ar@{->>}[d]\\
\Z/M\Z    \ar@{^{(}->}[r]    &\O_K/M\O_K.
}}

(iii) $\Longrightarrow$ (i). It suffice to show that under the hypothesis that every prime divisor $q$ of $M$ is such that $q>\sqrt M$. Let $q$ be a prime divisor of $M$ and $\Q$ a prime ideal of $\O_L$ lying over $q$ and let $\q=\Q\cap\O_K$. Clearly by the commutative diagram
\eqn{\xymatrix{
\Z_{(M)}    \ar@{->>}[d]\ar@{^{(}->}[r]    &\Z_{(q)}    \ar@{->>}[d]\\
\Z/M\Z    \ar@{->>}[r]    &\Z/q\Z
}}
we know that \eqref{MainThm.C3} also holds modulo $q$, i.e.,
\eqn{
T_{n-1}^{(m)}=\sum_{1\le i_1<\dots<i_m\le r}\prod_{j=1}^m t_{n-1}^{\s_{i_j}}\equiv\sum_{1\le i_1<\dots<i_m\le r}\prod_{j=1}^m u_{\xi}^{\s_{i_j}}=U^{(m)}\pmod q,\qquad m=1,\dots,r,
}
from which, noting the commutative diagram
\eqn{\xymatrix{
\Z_{(q)}    \ar@{->>}[d]\ar@{^{(}->}[r]    &(\O_K)_{\q}    \ar@{->>}[d]\\
\Z/q\Z    \ar@{^{(}->}[r]    &\O_K/\q,
}}
we know that
\eq{\label{EqUm}
z^r-U^{(1)}z^{r-1}+\dots+(-1)^rU^{(r)}=\prod_{i=1}^r(z-t_{n-1}^{\s_i})=\prod_{i=1}^r(z-u_{\xi}^{\s_i})
}
holds over the field $\O_K/\q$. Therefore $t_{n-1}=t_{n-1}^{\s_1}\equiv u_{\xi}^{\s_i}\pmod{\q}$ for some $i$ with $1\le i\le r$. Note that $t_{n-1}=\tau^{p^{n-1}}+\bar\tau^{p^{n-1}}$, $\tau\bar\tau=1$ and $u_{\xi}^{\s_i}=\xi^{\s_i}+\bar\xi^{\s_i}$. Still as \eqref{EqUm}, we obtain over $\O_L/\Q$ that
\eqn{
z^2-t_{n-1}z+1=(z-\tau^{p^{n-1}})(z-\bar\tau^{p^{n-1}})=(z-\xi^{\s_i})(z-\bar\xi^{\s_i})=(z-\xi^i)(z-\xi^{-i}),
}
hence $\tau^{p^{n-1}}\equiv\xi^{\pm i}\pmod{\Q}$ for some $i$ with $1\le i\le r$. That is, $\tau^{p^{n-1}}=\xi^{\pm i}$ when viewed in $\O_L/\Q$. So the order of $\tau$ in $(\O_L/\Q)^*$ is $p^n$. Consequently, $p^n\mid\card{(\O_L/\Q)^*}=\mathcal{N}\Q-1=q^{f(\Q|q)}-1\mid q^{p-1}-1$ since $f(\Q|q)\mid[L:\QQ]=p-1$. In other words, $q^{p-1}\equiv1\pmod{p^n}$. Since by hypothesis no solution of the last congruence equation greater than $1$ and less than $p^n$ is a divisor of $M$, it follows that $q\ge p^n>\sqrt M$. This completes the proof.
}

\section{Computational Considerations}\label{SecComp}
    Some care is needed here when applying the main theorem into computation.

\subsection{Solve an equation}\label{SecComp.Solve}
    Both in Proposition \ref{TrghoutProp} and Theorem \ref{MainThm} we have to solve the equation
\eq{\label{xp1modpn}
x^{p-1}\equiv1\pmod{p^n}.
}
Note that $(\Z/p^n\Z)^*$ is cyclic then the equation has exactly $p-1$ roots. In fact we can find a primitive root modulo $p^n$ easily (see \cite[Section 10.6]{AnalyticNT}). That is for a primitive root $g$ modulo $p$, if $g^{p-1}\not\equiv1\pmod{p^2}$ then $g$ is also a primitive root modulo $p^n$ for all $n\geq1$, otherwise $g+p$ makes it. Now suppose we find a primitive root $g$ modulo $p^n$, then all the $p-1$ roots of the equation \eqref{xp1modpn} are given by $g^{kp^{n-1}}\mod{p^n}$ for $k=0,\dots,p-2$. Set $w_{n}^{(k)}=g^{kp^{n-1}}\mod{p^n}$ with $0<w_{n}^{(k)}<p^n$ for $k=0,\dots,p-2$. Obviously, we have $w_{n}^{(0)}=1$ and $w_{n}^{\left(\frac{p-1}{2}\right)}=p^n-1$. Once we find such $g$, we can compute the values of $w_{n}^{(k)}$ for $k=1,\ldots,p-2$ and $k\neq\frac{p-1}{2}$.
Below we list the values of $g$ for $p<100$.

\begin{center}
\small{\textbf{Table 1} The values of primitive root $g$ modulo $p^n$ with $p<100$}
\begin{tabular}{|c|c|c|c|c|c|c|c|c|c|c|c|c|c|c|c|}
  \hline
  $p$ & $g$&$p$ & $g$ &$p$ & $g$ & $p$ & $g$&$p$ & $g$ &$p$ & $g$ &$p$ & $g$ &$p$ & $g$\\\hline
  3 & 2 & 11 & 2 & 19 & 2 & 31 & 3 & 43 & 3 & 59 & 2 & 71 & 7 & 83 & 2\\\hline
  5 & 2 & 13 & 2 & 23 & 5 & 37 & 2 & 47 & 5 & 61 & 2 & 73 & 5 & 89 & 3\\\hline
  7 & 3 & 17 & 3 & 29 & 2 & 41 & 6 & 53 & 2 & 67 & 2 & 79 & 3 & 97 & 5\\\hline
\end{tabular}
\end{center}

\subsection{Find $\pi$}\label{SecComp.Find}
    We now describe how to find the desired $\pi$. Suppose $p\leq19$ so it is well known that $\O_L=\Z[\zeta_p]$ is a Principal Ideal Domain(PID for short)(see \cite[Chapter 11]{Washington}). For $M$ as described before we find a small prime $l\equiv 1\pmod p$ such that $l\nmid M$ and $M$ is not a $p\th$ power modulo $l$, i.e. $M^{(l-1)/p}\not\equiv1\pmod{l}$ (by Extended Riemann Hypothesis, this could be found within $2(\log_2(M))^2$, assuming $M$ is not a $p^{\th}$ power, see \cite{Stein}). Let $\L$ be an ideal of $\O_L$ lying over $l$ and suppose $\L=\pi\O_L$. We can assume $\pi$ is primary (see Definition \ref{Primary}). Since $l\equiv 1\pmod p$, we have $f(\L|l)=1$ and then $\Z/l\Z\cong\O_L/\L$. It follows that $M$ is not a $p\th$ power modulo $\L$, i.e. $\fracp{M}{\pi}\ne1$ by Proposition \ref{pResSym1}(i). Since $\pi$ is primary we may use the Eisenstein reciprocity law (Theorem \ref{Eisenstein}) to obtain $\fracp{\pi}{M}=\fracp{M}{\pi}$ is a primitive $p\th$ root of unity. This gives an easy method to find $\pi$ when $p\leq19$. We will show the computational details in Section \ref{SecImple} latter. Here we give some examples of the values of $l$ and $\pi$ obtained during implementation (Section \ref{SecImple})(where $i$ stands for $w_n^{(i)}=g^{ip^{n-1}} \mod p^n$):
\begin{center}
\small{\textbf{Table 2} The values of $l$ and $\pi$ for $M=Ap^n+w_n^{(i)}$}
\begin{tabular}{|c|c|c|c|c|c|}
  \hline
  $A$ & $p$ & $n$ & $i$ & $l$ & $\pi$		\\\hline
  1&     3&    1&   1&  7& $1+3\zeta_3$\\\hline
  10 & 3 & 100 & 1 & 13 & $-4-3\zeta_3$\\\hline
  2 & 7 & 5 & 4 & 43 & $-1+\zeta_7-\zeta_7^3-\zeta_7^5 $	\\\hline
  10 & 7 & 100 & 4 & 29 & $-1-\zeta_7-2\zeta_7^2-\zeta_7^4-\zeta_7^5 $	\\\hline
\end{tabular}
\end{center}

    For $p\geq23$, we know that $\O_L=\Z[\zeta_p]$ is not necessarily a PID. Below we briefly describe a method to find $\pi$ due to P. Berrizbeitia et al. \cite{Amn+wn}, which contains details, see \cite{Amn+wn}. First, find a prime $q\equiv1\pmod{p}$ such that $u=M^{(q-1)/p}$ has order $p\pmod{q}$. Next, choose $\pi$ from the ideal $\Q\subseteq\O_L$ generated by $q$ and $\zeta_p-u$, compute the norm of $\pi$ and see if it satisfies the following condition: $\cnm{\pi}=tq$ with $p\nmid t$ and every prime divisor $l$ of $t$ satisfies $M^{(l-1)/\text{gcd}(l-1,p)}\equiv1\pmod{l}$. If so, then $\fracp{M}{\pi}$ is a primitive $p\th$ root of unity. Also, we may assume $\pi$ is primary, then use the Eisenstein reciprocity law (Theorem \ref{Eisenstein}) to obtain $\fracp{\pi}{M}=\fracp{M}{\pi}$ is a primitive $p\th$ root of unity.

    For some algorithms in algebraic number fields needed in the above methods, one can see \cite{Cohen}.

\subsection{Compute the values of $U^{(m)}$ for $m=1,\ldots,r$}\label{SecComp.Um}
    We could compute $U^{(m)}$ by definition \eqref{Um} directly. Another method is to compute the minimal polynomial of $\xi+\xi^{-1}$ over $\QQ$, by noting that \eqref{EqUm} implies $(-1)^mU^{(m)}$ is the coefficient of $z^{r-m}$ in the minimal polynomial. We give the list of $\{U^{(m)}\mid m=1,\ldots,r\}$ for $p\le19$ which is obtained during implementation (Section \ref{SecImple}).
\begin{center}
\small{\textbf{Table 3} The values of $U^{(m)}$ with $p\le19$}

\begin{tabular}{|c|l|}
  \hline
  $p$ & $\{U^{(m)}\mid m=1,\ldots,r\}$		\\\hline
  3   & $\{-1\}$				\\\hline
  5   & $\{-1,-1\}$				\\\hline
  7   & $\{-1,-2,1\}$				\\\hline
  11  & $\{-1,-4,3,3,-1\}$			\\\hline
  13  & $\{-1,-5,4,6,-3,-1\}$			\\\hline
  17  & $\{-1,-7,6,15,-10,-10,4,1\}$		\\\hline
  19  & $\{-1,-8,7,21,-15,-20,10,5,-1\}$	\\\hline
\end{tabular}
\end{center}

\subsection{Compute the polynomials $F(z_1,z_2)$ and $G_m(z_1,\dots,z_r)$ for $m=1,\ldots,r$}\label{SecComp.F.G}
    We remark here that these polynomials are done in pre-computation. The computation is standard symmetric polynomial reduction. With the help of \nolinkurl{SymmetricReduction[]} in \emph{Mathematica} \cite{Mathematica}, we obtain the following results.

    For $p=3$, $r=1$ the case is trivial: $F(z_1,z_2)=z_1^3-3 z_1 z_2$ and $G_1(z_1)=-3 z_1+z_1^3$.

    For $p=5$, we have $r=2$, $F(z_1,z_2)=z_1^5-5 z_1^3 z_2+5 z_1 z_2^2$ and
\aln{
G_1(z_1,z_2)&=5 z_1-5 z_1^3+z_1^5+15 z_1 z_2-5 z_1^3 z_2+5 z_1 z_2^2,\\
G_2(z_1,z_2)&=25 z_2-25 z_1^2 z_2+5 z_1^4 z_2+50 z_2^2-20
z_1^2 z_2^2+35 z_2^3-5 z_1^2 z_2^3+10 z_2^4+z_2^5.
}
This coincides with the results in P. Berrizbeitia et al. \cite{A5n}.

    For $p=7$, then $r=3$, $F(z_1,z_2)=z_1^7-7 z_1^5 z_2+14 z_1^3 z_2^2-7 z_1 z_2^3$ and

\aln{
G_1(z_1,z_2,z_3)&=-7 z_1+14 z_1^3-7 z_1^5+z_1^7-42 z_1 z_2+35 z_1^3 z_2 -7 z_1^5 z_2-35 z_1 z_2^2\\
	&+14 z_1^3 z_2^2-7 z_1 z_2^3+42 z_3-35 z_1^2 z_3+7 z_1^4 z_3+35 z_2 z_3-21 z_1^2 z_2 z_3+7 z_2^2 z_3\\
	&+7 z_1 z_3^2,}

\aln{
G_2(z_1,z_2,z_3)&=49 z_2-98 z_1^2 z_2+49 z_1^4 z_2-7 z_1^6 z_2+196 z_2^2-196 z_1^2 z_2^2+42 z_1^4 z_2^2+294 z_2^3\\
	&-161 z_1^2 z_2^3+14 z_1^4 z_2^3+210 z_2^4-56 z_1^2 z_2^4+77 z_2^5-7 z_1^2 z_2^5+14 z_2^6+z_2^7\\
	&+98 z_1 z_3-49 z_1^3 z_3+7 z_1^5 z_3-245 z_1 z_2 z_3+217 z_1^3 z_2 z_3-42 z_1^5 z_2 z_3\\
	&-469 z_1 z_2^2 z_3+168 z_1^3 z_2^2 z_3-273 z_1 z_2^3 z_3+35 z_1^3 z_2^3 z_3-70 z_1 z_2^4 z_3-7 z_1 z_2^5 z_3\\
	&+441 z_3^2-259 z_1^2 z_3^2+42 z_1^4 z_3^2+630 z_2 z_3^2-91 z_1^2 z_2 z_3^2-35 z_1^4 z_2 z_3^2+329 z_2^2 z_3^2\\
	&+35 z_1^2 z_2^2 z_3^2+77 z_2^3 z_3^2+14 z_1^2 z_2^3 z_3^2+7 z_2^4 z_3^2-91 z_1 z_3^3+35 z_1^3 z_3^3-91 z_1 z_2 z_3^3\\
	&-7 z_1^3 z_2 z_3^3-21 z_1 z_2^2 z_3^3+21 z_3^4+7 z_1^2 z_3^4+7 z_2 z_3^4,}

\aln{
G_3(z_1,z_2,z_3)&=-343 z_3+686 z_1^2 z_3-343 z_1^4 z_3+49 z_1^6 z_3-1372 z_2 z_3+1372 z_1^2 z_2 z_3\\
	&-294 z_1^4 z_2 z_3-2058 z_2^2 z_3+1127 z_1^2 z_2^2 z_3-98 z_1^4 z_2^2 z_3-1470 z_2^3 z_3+392 z_1^2 z_2^3 z_3\\
	&-539 z_2^4 z_3+49 z_1^2 z_2^4 z_3-98 z_2^5 z_3-7 z_2^6 z_3+1372z_1 z_3^2-1078 z_1^3 z_3^2+196 z_1^5 z_3^2\\
	&+2156 z_1 z_2 z_3^2-784 z_1^3 z_2 z_3^2+1372 z_1 z_2^2 z_3^2-196 z_1^3 z_2^2 z_3^2+392 z_1 z_2^3 z_3^2\\
	&+42z_1 z_2^4 z_3^2+833 z_3^3-1176 z_1^2 z_3^3+294 z_1^4 z_3^3+1176 z_2 z_3^3-784 z_1^2 z_2 z_3^3+637 z_2^2 z_3^3\\
	&-161 z_1^2 z_2^2 z_3^3+154 z_2^3 z_3^3+14 z_2^4 z_3^3-490 z_1 z_3^4+210 z_1^3 z_3^4-308 z_1 z_2 z_3^4\\
	&-56 z_1 z_2^2 z_3^4-70 z_3^5+77 z_1^2 z_3^5-42 z_2 z_3^5-7 z_2^2 z_3^5+14 z_1 z_3^6+z_3^7.
}

We will see more details of this case in Section \ref{Sec.p=7}.

    For $p\ge11$, $G_m$ is too long to be written down here. We just give a list of $F(z_1,z_2)$.
\begin{center}
\small{\textbf{Table 4} The values of $F(z_1,z_2)$}
\begin{tabular}{|c|l|}
  \hline
  $p$ & $F(z_1,z_2)$\\\hline
  11  & $z_1^{11}-11 z_1^9 z_2+44 z_1^7 z_2^2-77 z_1^5 z_2^3+55 z_1^3 z_2^4-11 z_1 z_2^5$\\\hline
  13  & $z_1^{13}-13 z_1^{11} z_2+65 z_1^9 z_2^2-156 z_1^7 z_2^3+182 z_1^5 z_2^4-91 z_1^3 z_2^5+13 z_1 z_2^6$\\\hline
  17  & $z_1^{17}-17 z_1^{15} z_2+119 z_1^{13} z_2^2-442 z_1^{11} z_2^3+935 z_1^9 z_2^4-1122 z_1^7 z_2^5+714 z_1^5 z_2^6-204 z_1^3 z_2^7+17 z_1 z_2^8$\\\hline
  19  & $z_1^{19}-19 z_1^{17} z_2+152 z_1^{15} z_2^2-665 z_1^{13} z_2^3+1729 z_1^{11} z_2^4-2717 z_1^9 z_2^5+2508 z_1^7 z_2^6-1254 z_1^5 z_2^7$\\
      & $+285 z_1^3 z_2^8-19 z_1 z_2^9$\\\hline
\end{tabular}
\end{center}

\section{A Concrete Example When $p=7$}\label{Sec.p=7}
    For $p=7$ we rewrite the results achieved before so that one can see the applications in computation. This was done by P. Berrizbeitia et al. \cite{A5n} when $p=5$.

    Suppose $p=7$, $M=A7^n+w_n$ and we do not change other notations. Thus $f=\ord{7}{M}\mid [L:\QQ]=6$ so we know that the only possible values of $f$ are $1$, $2$, $3$ and $6$, and correspondingly, select $\gamma=1-\s_{-1}+2(\s_4-\s_3)+3(\s_5-\s_2$), $1+2\s_4+3\s_5$, $1-\s_{-1}$ and $1$ (see \eqref{gamma1} and \eqref{gamma2}). By writing each $\Phi_f(x)$ out, we know $\tau=\pi^{\gamma\frac{M-1}{7^n}}$, $\pi^{\gamma(\s_{M}-1)\frac{M+1}{7^n}}$, $\pi^{\gamma(\s_M-1)\frac{M^2+M+1}{7^n}}$ and $\pi^{\gamma(\s_{M}-1)(\s_M+1)(\s_M^2+\s_M+1)\frac{M^2-M+1}{7^n}}$ respectively (see \eqref{tau}). Here $\pi$ could be found by the method described in Section \ref{SecComp.Find}. That is, find a small prime $l\equiv 1\pmod 7$ such that $M$ is not a $7\th$ power modulo $l$. Let $\L$ be an ideal of $\O_L$ lying over $l$ and find a primary $\pi$ such that $\L=\pi\O_L$ (since $\O_L$ is a PID). Since $r=3$ we write $T_k$, $J_k$, and $N_k$ for $T_k^{(m)},\,m=1,2,3$ and we see \eqref{Tkm} in details that
\aln{
T_k&=\tr{K}{\QQ}{t_k}=\sum_{i=1}^3t_k^{\s_i}, \\
J_k&=\nm{K}{\QQ}{t_k}\tr{K}{\QQ}{t_k^{-1}}=\sum_{1\le i<j\le3}t_k^{\s_i}t_k^{\s_j}, \\
N_k&=\nm{K}{\QQ}{t_k}=\prod_{i=1}^3t_k^{\s_i}.
}
also $U^{(1)}=-1$, $U^{(2)}=-2$ and $U^{(3)}=1$.
Now the main theorem becomes
\thm{
Let $M$, $\tau$, $t_k$, $T_k$, $J_k$ and $N_k$ be as before. Suppose further that $M$ is not divisible by any of the solutions of $x^6\equiv1\pmod{7^n}$; $1<x<7^n$. Then the following statements are equivalent:
\begin{enumerate}[(i)]
\item $M$ is prime.
\item There exists a primitive $7\th$ root of unity $\xi$ such that
\eqn{
\tau^{7^{n-1}}\equiv\xi\pmod M.
}
\item
\aln{
T_k&\equiv-1\pmod M,\\
J_k&\equiv-2\pmod M,\\
N_k&\equiv1\pmod M.
}
\end{enumerate}
}
With the computational results obtained in Sections \ref{SecComp.Um} and \ref{SecComp.F.G}, we give the recursive formulas
\aln{
T_0&=\tr{K}{\QQ}{\tau+\bar\tau}, \\
J_0&=\nm{K}{\QQ}{\tau+\bar\tau}\tr{K}{\QQ}{(\tau+\bar\tau)^{-1}},\\
N_0&=\nm{K}{\QQ}{\tau+\bar\tau}.
}
And for $k\ge0$,
\aln{\sp{
T_{k+1} &= 42 N_k+35 J_k N_k+7 J_k^2 N_k-7 T_k-42 J_k T_k-35 J_k^2 T_k-7 J_k^3 T_k+7 N_k^2 T_k-35 N_k T_k^2\\
	&-21 J_k N_k T_k^2+14 T_k^3+35 J_k T_k^3+14	J_k^2 T_k^3+7 N_k T_k^4-7 T_k^5-7 J_k T_k^5+T_k^7,
}}
\aln{\sp{					
J_{k+1} &= 49 J_k+196 J_k^2+294 J_k^3+210 J_k^4+77 J_k^5+14 J_k^6+J_k^7+441 N_k^2+630 J_k N_k^2\\
	&+329 J_k^2 N_k^2+77 J_k^3 N_k^2+7 J_k^4 N_k^2+21 N_k^4+7 J_k N_k^4+98 N_k T_k-245 J_k N_k T_k\\
	&-469 J_k^2 N_k T_k-273 J_k^3 N_k T_k-70 J_k^4 N_k T_k-7 J_k^5 N_k T_k-91 N_k^3 T_k-91 J_k N_k^3 T_k\\
	&-21 J_k^2 N_k^3 T_k-98 J_k T_k^2-196 J_k^2 T_k^2-161 J_k^3 T_k^2-56 J_k^4 T_k^2-7 J_k^5 T_k^2\\
	&-259 N_k^2 T_k^2-91 J_k N_k^2 T_k^2+35 J_k^2 N_k^2 T_k^2+14 J_k^3 N_k^2 T_k^2+7 N_k^4 T_k^2\\
	&-49 N_k T_k^3+217 J_k N_k T_k^3+168 J_k^2 N_k T_k^3+35 J_k^3 N_k T_k^3+35 N_k^3 T_k^3\\
	&-7 J_k N_k^3 T_k^3+49 J_k T_k^4+42 J_k^2 T_k^4+14 J_k^3 T_k^4+42 N_k^2 T_k^4-35 J_k N_k^2 T_k^4\\
	&+7 N_k T_k^5-42 J_k N_k T_k^5-7 J_k T_k^6,
}}
\aln{\sp{	
N_{k+1} &= -343	N_k-1372 J_k N_k-2058 J_k^2 N_k-1470 J_k^3 N_k-539 J_k^4 N_k-98 J_k^5 N_k-7 J_k^6 N_k\\
	&+833 N_k^3+1176 J_k N_k^3+637 J_k^2 N_k^3+154 J_k^3 N_k^3+14 J_k^4 N_k^3-70 N_k^5-42 J_k N_k^5\\
	&-7 J_k^2 N_k^5+N_k^7+1372 N_k^2 T_k+2156 J_k N_k^2 T_k+1372 J_k^2 N_k^2 T_k+392 J_k^3 N_k^2 T_k\\
	&+42 J_k^4 N_k^2 T_k-490 N_k^4 T_k-308 J_k N_k^4 T_k-56 J_k^2 N_k^4 T_k+14 N_k^6 T_k+686 N_k T_k^2\\
	&+1372 J_k N_k T_k^2+1127 J_k^2 N_k T_k^2+392 J_k^3 N_k T_k^2+49 J_k^4 N_k T_k^2-1176 N_k^3 T_k^2\\
	&-784 J_k N_k^3 T_k^2-161 J_k^2 N_k^3 T_k^2+77 N_k^5 T_k^2-1078 N_k^2 T_k^3-784 J_k N_k^2 T_k^3\\
	&-196 J_k^2 N_k^2 T_k^3+210 N_k^4 T_k^3-343 N_k T_k^4-294 J_k N_k T_k^4-98 J_k^2 N_k T_k^4\\
	&+294 N_k^3 T_k^4+196 N_k^2 T_k^5+49 N_k T_k^6.
}}
Besides, by Section \ref{SecComp.Solve} we are able to choose $3$ as a primitive root modulo $7^n$ to solve \eqref{xp1modpn}.

\section{Implementation and Computational results}\label{SecImple}
    We use \emph{PARI} \cite{PARI}, a widely used computer algebra system designed for fast computation in number theory originally developed by Henri Cohen and his co-workers, to implement the algorithm in C language. We run our program on a computer with Intel Xeon E5530 2.40GHz CPU and 96GB memory.

	For given $M$ with $p\le19$, there is no difficult to compute required $l$, by searching primes less than $2(\log_2(M))^2$ that is congruent to $1$ modulo $p$ as described in Section \ref{SecComp.Find}. We remark that such an $l$ was found for all of our tested $M$'s except for some $M<100$. Compute $U^{(m)}$ by the minimal polynomial of $\xi+\xi^{-1}$ over $\QQ$. By prime decomposition we can find a generator of the ideal $\pi\O_L$, for the primary associate of which we use method described in \cite[Section 14.2]{GTM84}. To compute $\tau$, we note that it is no harm to reduce the intermediate computation result modulo $M$ (this may make the failure of Lemma \ref{Neq1} but the lemma is only used to deduce the recursive formulas), thus we could launch the fast power modulo $M$ to obtain $\pi^{\Phi_f(M)/p^n}$ in at most $2p\log_2(M)$ multiplications of $\pi$ modulo $M$ and the latter could be viewed as integer arithmetic modulo $M$ since there is a standard integral basis for $\O_L$. This could be omitted when $w=\pm1$. The power to $\gamma\prod_{d\mid f,d<f}\Phi_d(\s_M)\th$ is easy since $p$ is small. The polynomials $G_m(x_1,\dots,x_r)$ are obtained as pre-computation (see Section \ref{SecComp.F.G}). By \eqref{Tkm} we find that $(-1)^mT_0^{(m)}$ is the coefficient of $z^{r-m}$ in the characteristic polynomial of $t_0$ over $K/\QQ$ and the subsequent $T_k^{(m)}$ is obtained by integer arithmetic modulo $M$ using the recursive formulas in Proposition \ref{PropRecurT} for $n-1$ ($<\log_p(M)$) times. The final work is to verify the congruences \eqref{MainThm.C3} and to make a further check of $p-1$ solutions of \eqref{xp1modpn}. The total complexity is a polynomial of $\log_2(M)$ if $p$ is viewed as a constant. The computational results are briefly described as follows.

    We first verified the correctness (i.e. the primarily of $M$ given by our program) of the algorithm in the small range $p=3$, $5$, $0\le A\le100$, $1\le n\le1000$ and all $p-1$ values of $i$, finding no mistakes. We also verified all numbers for $p=7$ in the range $0\le A\le8$ and $1\le n\le1000$. There are $3263$ primes of all the verified $658732$ numbers. Some of them along with the cost times are listed in the following table:
\begin{center}
\small{\textbf{Table 5} The running time of some verified numbers}
\begin{tabular}{|c|c|c|c|c|l|l|}
  \hline
  $A$ & $p$ & $n$& $i$ & Primality & Time (ms) \\\hline
  1&3&1&1&yes&3 \\\hline
  100&3&911&0&yes&8 \\\hline
  100&3&1000&1&no&11 \\\hline
  2&5&100&0&no&9 \\\hline
  3&5&171&2&yes&11 \\\hline
  3&5&1000&3&no&437 \\\hline
  100&5&992&3&yes&436 \\\hline
  0&7&1&1&yes&39 \\\hline
  3&7&984&4&no&2732 \\\hline
  8&7&806&1&yes&1540 \\\hline
  8&7&1000&5&no&2538 \\\hline
\end{tabular}
\end{center}
Among all numbers in this range, $M=3\times7^{984}+w_{984}^{(4)}$ takes the longest time 2732 ms. And we also noted that for large primes the APRCL test took much more time than our algorithm. For instance, the prime $M=8\times7^{806}+w_{806}^{(1)}$ takes $1540$ms using our algorithm but $140892$ms using APRCL test. Next we used our program to search all $6060$ numbers $M$ in the range $p=7$, $1\le A\le10$, $2000\le n\le2100$ and all $i$, finding the only two primes are $7\times7^{2077}+w_{2077}^{(5)}$ and $8\times7^{2060}+w_{2060}^{(5)}$. The longest time of these numbers taken are 17472ms. Another test was aimed to the only two BPSW-pseudoprimes (see \cite{PSW}) $M$ in the range $p=7$, $A=1$, $5000\le n\le7000$ and all $i$, i.e. $7^{5180}+w_{5180}^{(3)}$ and $7^{5618}+w_{5618}^{(2)}$. Our program asserted that they are all primes, taking $166786$ms and $319407$ms respectively. However we terminated the program using APRCL test to process these two large primes because we had waited for more than three days.

\vspace{5mm}

\noindent \textbf{Acknowledgments}\quad  The work of this paper
was supported by the NNSF of China (Grants Nos. 11071285,
61121062), 973 Project (2011CB302401) and the National Center for Mathematics and Interdisciplinary Sciences, CAS.

\end{document}